\documentclass{amsproc}
\usepackage[dvips]{graphicx}
\usepackage{amsmath,amssymb}

\begin{document}

\title[On an inverse dynamic problem for the wave equation.]
{Inverse dynamic problem for the wave equation with periodic
boundary conditions.}

\author[A.\,S.~Mikhaylov, V.\,S.~Mikhaylov]
{$^1$, $^2$A.\,S.~Mikhaylov, $^1$, $^2$V.\,S.~Mikhaylov}

\address{
$^1$ St. Petersburg Department of V.A. Steklov Institute of
Mathematics of the Russian Academy of Sciences, 7, Fontanka,
191023 St. Petersburg, Russia. $^2$ Saint Petersburg State
University, St.Petersburg State University, 7/9 Universitetskaya
nab., St. Petersburg, 199034 Russia.}

\email{mikhaylov@pdmi.ras.ru, vsmikhaylov@pdmi.ras.ru}

\begin{abstract}
We consider the inverse dynamic problem for the wave equation with
a potential on an interval $(0,2\pi)$ with periodic boundary
conditions. We use a boundary triplet to set up the
initial-boundary value problem. As an inverse data we use a
response operator (dynamic Dirichlet-to-Neumann map). Using the
auxiliary problem on the whole line, we derive equations of the
inverse problem. We also establish the relationships between
dynamic and spectral inverse data.
\end{abstract}


\maketitle

\newtheorem{corollary}{Corollary}
\newtheorem{definition}{Definition}
\newtheorem{lemma}{Lemma}
\newtheorem{proposition}{Proposition}
\newtheorem{remark}{Remark}
\newtheorem{theorem}{Theorem}

\section{Introduction}

Inverse problems for one-dimensional continuous and discrete
systems plays an important role for the creation of new
nano-devices, to mention just \cite{R,V} and references therein.
In the present paper we set up and study the inverse dynamic
problem for a wave equation with a potential on an interval with a
periodic boundary conditions. The control problems for dynamical
systems for wave equation with periodic boundary conditions (the
density allows certain dependence on time) were considered in
\cite{ABI,ABP}. The spectral problem for a Schr\"odinger operator
on an interval with periodic and anti-periodic boundary conditions
are used for treating the spectral problem for a Schr\"odinger
operator with periodic potential on $\mathbb{R}$, see \cite{Le}.
The inverse spectral problem with periodic boundary conditions for
Schr\"odinger operator plays an important role for studying
inverse problems on graphs with cycles \cite{K}.

In the previous papers by the authors the "dynamic" approach to
inverse spectral problems based on ideas of the Boundary Control
method \cite{B07,B17} was developed in the cases of Schr\"odinger
operator on a half-line \cite{AMR,MM3,MM6,MM7} and finite and
semi-infinite Jacobi matrices \cite{MMS,MMS,MM8}. We believe that
our "dynamic" methods will help us to establish new relationships
and develop new tools for studying the inverse problems with
periodic potential, and will also stimulate studying inverse
problems on graphs with cycles \cite{BW,K}.

For a potential $q\in C^2(0,2\pi)$ we consider an operator $H$ in
$L_2(0,2\pi)$ given by
\begin{eqnarray*}
(Hf)(x)=-f''(x)+q(x)f(x),\quad x\in (0,2\pi),\\
\operatorname{dom}H=\left\{f\in H^2(0,2\pi)\,|\,
f(0)=f'(0)=f(2\pi)=f'(2\pi)=0\right\}.
\end{eqnarray*}
Then
\begin{eqnarray*}
(H^*f)(x)=-f''(x)+q(x)f(x),\quad x\in (0,2\pi),\\
\operatorname{dom}H^*=\left\{f\in H^2(0,2\pi)\right\}.
\end{eqnarray*}
For a continuous function $g$ we introduce the notations
\begin{equation*}
g_{0}:=\lim_{\varepsilon\to 0}g(0+\varepsilon),\quad
g_{2\pi}:=\lim_{\varepsilon\to 0}g(2\pi-\varepsilon).
\end{equation*}
Let $B:=\mathbb{R}^2$. The \emph{boundary operators}
$\Gamma_{0,1}: \operatorname{dom}H^*\mapsto B$ are introduced by
the rules
\begin{equation*}
\Gamma_0w:=\begin{pmatrix}w_0-w_{2\pi}\\
w'_0-w'_{2\pi}\end{pmatrix},\quad
\Gamma_1w:=\frac{1}{2}\begin{pmatrix}w'_{0}+w'_{2\pi}\\
-w_{0}-w_{2\pi}\end{pmatrix}.
\end{equation*}
Integrating by parts for $u,v\in \operatorname{dom}H^*$ shows that
the abstract second Green identity holds:
\begin{equation*}
\left(H^*u,v\right)_{L_2(0,2\pi)}-\left(u,H^*v\right)_{L_2(0,2\pi)}=\left(\Gamma_1u,\Gamma_0v\right)_B-\left(\Gamma_0u,\Gamma_1v\right)_B.
\end{equation*}
The mapping
\begin{equation*}
\Gamma:=\begin{pmatrix}\Gamma_0\\\Gamma_1\end{pmatrix}:\operatorname{dom}H^*\mapsto
B\times B
\end{equation*}
is surjective. Then a triplet $\{B,\Gamma_0,\Gamma_1\}$ is a
\emph{boundary triplet} for $H^*$ (see \cite{BMN}).

Let $T>0$ be fixed. We use the triplet $\{B,\Gamma_0,\Gamma_1\}$
to set up the initial-boundary value problem for dynamical system:
\begin{equation}
\label{wave_eqn}
\begin{cases}
u_{tt}+H^*u=0,\quad t> 0,\\
(\Gamma_0u)(t)=\begin{pmatrix}f_1(t) \\ f_2'(t)\end{pmatrix}, \quad t>0,\\
u(\cdot,0)=u_t(\cdot,0)=0.
\end{cases}
\end{equation}
Here the vector function $F=\begin{pmatrix}f_1 \\
f_2\end{pmatrix}$, $f_1\, f_2\in L_2(0,T),$ is interpreted as a
\emph{boundary control}. The solution to (\ref{wave_eqn}) is
denoted by $u^F.$ The \emph{response operator} is introduced by
the rule
\begin{equation*}
\left(R^TF\right)(t):=\left(\Gamma_1u^F\right)(t),\quad t>0.
\end{equation*}
The speed of the wave propagation in the system (\ref{wave_eqn})
equal to one, that is why the natural set up of the dynamic
inverse problem (IP) is to find a potential $q(x),$ $x\in
(0,2\pi)$ from the knowledge of a response operator $R^{2\pi}$
(cf. \cite{B07,B17,BM01,AM}).

In the second section we derive the representation formula for the
solution $u^F$, introduce the auxiliary dynamical system on the
real line (see also \cite{MM5}), and use the finiteness of the
speed of wave propagation to establish relationships between the
problem with periodic boundary conditions and problem on
$\mathbb{R}$. In the third section, on the basis of this
relationship, we obtain the suitable version of Krein and
Gelfand-Levitan equations of the dynamic inverse problem. In the
last section we derive the spectral representation of the response
operator and dynamic representation of a Weyl function associated
with $\{B,\Gamma_0,\Gamma_1\}$.

\section{Forward problem, auxiliary dynamical system.}

We introduce the \emph{outer space} of the system
(\ref{wave_eqn}), the space of controls as $\mathcal{F}^T:=
L_2(0,T;\mathbb{R}^2)$, $F\in \mathcal{F}^T$,
$F=\begin{pmatrix}f_1\\f_2\end{pmatrix}$. By $q$ we also denote
the same potential, periodically continued to the whole real line:
$q(x+2\pi)=q(x),$ $x\in \mathbb{R}$.
\begin{theorem}
The solution to (\ref{wave_eqn}) with a control $F\in
\mathcal{F}^T\cap C_0^\infty(\mathbb{R}_+)$, admits the following
representation:
\begin{itemize}
\item[1)] For $0<t<2\pi$
\begin{eqnarray}
\label{UF_repr} u^F(x,t)=u^{F_+}_1(x,t)+u^{F_-}_1(x,t)\\
=\frac{1}{2}f_1(t-x)-\frac{1}{2}f_2(t-x)
+\int_x^t w_1^0(x,s)f_1(t-s)+w_2^0(x,s)f_2(t-s)\,ds\notag\\
-\frac{1}{2}f_1(t+x-2\pi)-\frac{1}{2}f_2(t+x-2\pi)+\int_{2\pi-x}^t
w_1^{2\pi}(x,s)f_1(t-s)+w_2^{2\pi}(x,s)f_2(t-s)\,ds.\notag
\end{eqnarray}
where kernels $w_{1,2}^{0,2\pi}(x,t)$ satisfy the following
Goursat problems:
\begin{eqnarray}
\begin{cases}
{w_{1}^0}_{tt}(x,t)-{w_{1}^0}_{xx}(x,t)+q(x){w_{1}^0}(x,t)=0,\quad 0<x<t,\\
\frac{d}{dx}{w_{1}^0}(x,x)=-\frac{q(x)}{4},\quad x>0,\\
{w_{1}^{2\pi}}_{tt}(x,t)-{w_{1}^{2\pi}}_{xx}(x,t)+q(x){w_{1}^{2\pi}}(x,t)=0,\quad
0<2\pi-x<t,\\
\frac{d}{dx}{w_{1}^{2\pi}}(x,2\pi-x)=-\frac{q(x)}{4},\quad x>0,\\
w^0_1(0,s)=w^{2\pi}_1(2\pi,s),\\
{w^0_1}_x(0,s)={w^{2\pi}_1}_x(2\pi,s).
\end{cases}
\label{W1} \\
\begin{cases}
{w_{2}^0}_{tt}(x,t)-{w_{2}^0}_{xx}(x,t)+q(x){w_{2}^0}(x,t)=0,\quad 0<x<t,\\
\frac{d}{dx}{w_{2}^0}(x,x)=\frac{q(x)}{4},\quad x>0,\\
{w_{2}^{2\pi}}_{tt}(x,t)-{w_{2}^{2\pi}}_{xx}(x,t)+q(x){w_{2}^{2\pi}}(x,t)=0,\quad
0<2\pi-x<t,\\
\frac{d}{dx}{w_{2}^{2\pi}}(x,2\pi-x)=-\frac{q(x)}{4},\quad x>0,\\
w^0_2(0,s)=w^{2\pi}_2(2\pi,s),\\
{w^0_2}_x(0,s)={w^{2\pi}_2}_x(2\pi,s).
\end{cases} \label{W2}
\end{eqnarray}

\item[2)] On $0<t<4\pi$
\begin{eqnarray*}
u^F(x,t)=u^{F_+}_1(x,t)+u^{F_-}_1(x,t)+u^{F_+}_2(x,t)+u^{F_-}_2(x,t)\\
=\frac{1}{2}f_1(t-x)-\frac{1}{2}f_2(t-x)
+\int_x^t w_1^0(x,s)f_1(t-s)+w_2^0(x,s)f_2(t-s)\,ds\notag\\
-\frac{1}{2}f_1(t+x-2\pi)-\frac{1}{2}f_2(t+x-2\pi)+\int_{2\pi-x}^t
w_1^{2\pi}(x,s)f_1(t-s)+w_2^{2\pi}(x,s)f_2(t-s)\,ds\notag \\
+\frac{1}{2}f_1(t-2\pi-x)-\frac{1}{2}f_2(t-2\pi-x)\notag\\
+\int_x^{t-2\pi} \widetilde w_1^0(x,s)f_1(t-2\pi-s)+\widetilde w_2^0(x,s)f_2(t-2\pi-s)\,ds\notag\\
-\frac{1}{2}f_1(t+x-4\pi)-\frac{1}{2}f_2(t+x-4\pi)\notag\\
+\int_{2\pi-x}^{t-2\pi}
\widetilde w_1^{2\pi}(x,s)f_1(t-2\pi-s)+\widetilde
w_2^{2\pi}(x,s)f_2(t-2\pi-s)\,ds.\notag
\end{eqnarray*}
where the integral kernels $w^{0,2\pi}_{1,2},$ $\widetilde
w^{0,2\pi}_{1,2}$  satisfy certain Goursat problems and the
following compatibility conditions:
\begin{eqnarray*}
w^0_{1,2}(0,s)=w^{2\pi}_{1,2}(2\pi,s),\,\, {w^0_{1,2}}_x(0,s)={w^{2\pi}_{1,2}}_x(2\pi,s),\quad 0<s<4\pi,\\
w^0_{1,2}(2\pi,s)=\widetilde w^0_{1,2}(0,s-2\pi),\,\,{w^0_{1,2}}_x(2\pi,s)=\widetilde {w^0_{1,2}}_x(0,s-2\pi),\quad 0<s<4\pi,\\
w^{2\pi}_{1,2}(0,s)=\widetilde
w^{2\pi}_{1,2}(2\pi,s-2\pi),\,\,{w^{2\pi}_{1,2}}_x(0,s)=\widetilde
{w^{2\pi}_{1,2}}_x(2\pi,s-2\pi),\quad 0<s<4\pi.
\end{eqnarray*}

\item[3)] On $0<t<2n\pi$, $n>1:$
\begin{equation}
\label{Repr_3}
u^F(x,t)=u^{F_+}_1(x,t)+u^{F_-}_1(x,t)+\ldots+u^{F_+}_n(x,t)+u^{F_-}_n(x,t),
\end{equation}
where
\begin{eqnarray*}
u^{F_+}_k(x,t)=\frac{1}{2}f_1(t-x-2(k-1)\pi)-\frac{1}{2}f_2(t-x-2(k-1)\pi)\\
+\int_{x+2(k-1)\pi}^t w_1(x+2(k-1)\pi,s)f_1(t-s)+w_2(x+2(k-1)\pi,s)f_2(t-s)\,ds\\
u^{F_-}_k(x,t)=-\frac{1}{2}f_1(t+x-2k\pi)-\frac{1}{2}f_2(t+x-2k\pi)\\
+\int_{2k\pi-x}^t
w_1(x-2k\pi,s)f_1(t-s)+w_2(x-2k\pi,s)f_2(t-s)\,ds
\end{eqnarray*}
and kernels $w_{1,2}$ satisfy the following Goursat problem:
\begin{eqnarray}
\begin{cases}
{w_{1}}_{tt}(x,t)-{w_{1}}_{xx}(x,t)+q(x){w_{1}}(x,t),\quad 0<|x|<t<2n\pi,\\
\frac{d}{dx}{w_{1}}(x,x)=-\frac{q(x)}{4},\quad x>0,\\
\frac{d}{dx}{w_{1}}(x,-x)=-\frac{q(x)}{4},\quad x<0,
\end{cases}
\label{G1}\\
\begin{cases}
{w_{2}}_{tt}(x,t)-{w_{2}}_{xx}(x,t)+q(x){w_{2}}(x,t),\quad 0<|x|<t<2n\pi,\\
\frac{d}{dx}{w_{2}}(x,x)=\frac{q(x)}{4},\quad x>0,\\
\frac{d}{dx}{w_{2}}(x,-x)=-\frac{q(x)}{4},\quad x<0.
\end{cases} \label{G2}
\end{eqnarray}
\end{itemize}
\end{theorem}
Several remarks have to be made.
\begin{remark}
The proof of the representation (\ref{UF_repr}) is straightforward
and similar to one in \cite{MM5}. If $F\in \mathcal{F}^T$, the
function $u^F$ defined by (\ref{UF_repr}) is a generalized
solution to (\ref{wave_eqn}) for $t\in(0,2\pi)$.
\end{remark}

\begin{remark}
The compatibility conditions in (\ref{W1}), (\ref{W2}) is used in
the next subsection to relate the solution of the problem with
periodic boundary conditions with one of the problem on the whole
line.
\end{remark}
Since we consider the periodic boundary conditions, sometimes it
would be convenient for us to interpret the interval as a ring:
\begin{remark}
The compatibility conditions in 2) allows one to construct the
"general" Goursat problems in 3). The physical meaning of the
representation (\ref{Repr_3}) is clear: the members of the sum
indexed with $"+"$ corresponds to waves that move clockwise, ones
indexed with $"-"$ correspond to waves moving counterclockwise.
\end{remark}

The \emph{response operator} $R^T: \mathcal{F}^T\mapsto
\mathcal{F}^T$ with the domain $D_R=\left\{\mathcal{F}^T\cap
C_0^\infty(0,T;\mathbb{R}^2)\right\}$ is defined by the rule
\begin{equation*}
(R^TF)(t):=\left(\Gamma_1u^F\right)(t),\quad 0<t<T.
\end{equation*}
Representation (\ref{UF_repr}) implies the following
\begin{corollary}
The response operator has a form:
\begin{itemize}
\item[1)] on an interval $(0,2\pi)$:
\begin{eqnarray}
\label{Resp_repr} \left(R^TF\right)(t)=
-\frac{1}{2}\begin{pmatrix}f_1'(t)\\-f_2(t)\end{pmatrix}+R*\begin{pmatrix}f_1\\f_2\end{pmatrix}.
\end{eqnarray}
where
\begin{equation*}
R(t):=\begin{pmatrix}r_{11}(t) & r_{12}(t)\\
r_{21}(t) & r_{22}(t)\end{pmatrix}=\begin{pmatrix}{w_1^0}_x(0,t) & {w_2^0}_x(0,t)\\
-{w_1^0}(0,t) & -{w_2^0}(0,t)\end{pmatrix}=\begin{pmatrix}{w_1^{2\pi}}_x(0,t) & {w_2^{2\pi}}_x(0,t)\\
-{w_1^{2\pi}}(0,t) & -{w_2^{2\pi}}(0,t)\end{pmatrix}
\end{equation*}
is a \emph{response matrix}.

\item[2)] on an interval $(0,2n\pi):$
\begin{eqnarray}
\label{Resp_repr1} \left(R^TF\right)(t)=
\left(-\frac{1}{2}\sum_{k=1}^n\begin{pmatrix}
\delta'(t-2k\pi)&0\\0&-\delta(t-2k\pi)\end{pmatrix}+\widetilde
R(t)\right)*\begin{pmatrix}f_1\\f_2\end{pmatrix},
\end{eqnarray}
where the integral kernel $\widetilde R$ is expressed in terms of
solutions to Goursat problems (\ref{G1}), (\ref{G2}).
\end{itemize}
\end{corollary}
\begin{remark}
Due to the finite speed of wave propagation in system
(\ref{wave_eqn}), the natural set up of IP is to recover the
potential on $(0,2\pi)$ from $R^{2\pi}$, that is why, for solving
IP we can consider the system for times less or equal $2\pi$.
\end{remark}

\subsection{Auxiliary problem on $\mathbb{R}$}

We introduce the the potential $\widetilde q$ by the rule
\begin{equation}
\widetilde q(x)=\begin{cases} q(x),\quad 0<x<2\pi,\\
0,\quad x> 2\pi,\\
 q(x+2\pi),\quad -2\pi<x<0,\\
 0,\quad x<- 2\pi,\\
\end{cases}
\end{equation}
For this potential we consider an operator $\widetilde H$ in
$L_2(\mathbb{R})$ given by
\begin{eqnarray*}
(\widetilde Hf)(x)=-f''(x)+\widetilde q(x)f(x),\quad x\in \mathbb{R},\\
\operatorname{dom}\widetilde H=\left\{f\in H^2(\mathbb{R})\,|\,
f(0)=f'(0)=0\right\}.
\end{eqnarray*}
Then
\begin{eqnarray*}
({\widetilde H}^*f)(x)=-f''(x)+\widetilde q(x)f(x),\quad x\in \mathbb{R},\\
\operatorname{dom}{\widetilde H}^*=\left\{f\in
L_2(\mathbb{R})\,|\, f\in H^2(-\infty,0),\,f\in
H^2(-\infty,0)\right\}.
\end{eqnarray*}
For a continuous function $g$ we denote
\begin{equation*}
g_\pm:=\lim_{\varepsilon\to 0}g(0\pm\varepsilon).
\end{equation*}
The \emph{boundary operators} $\widetilde \Gamma_{0,1}:
\operatorname{dom}H^*\mapsto B$ are introduced by the rules
\begin{equation*}
\widetilde\Gamma_0w:=\begin{pmatrix}w_+-w_-\\
w'_+-w'_-\end{pmatrix},\quad
\widetilde\Gamma_1w:=\frac{1}{2}\begin{pmatrix}w'_++w'_-\\
-w_+-w'_-\end{pmatrix}.
\end{equation*}
We consider the initial boundary value problem for an auxiliary
dynamical system on $\mathbb{R}$:
\begin{equation}
\label{wave_eqn1}
\begin{cases}
v_{tt}+v_{xx}+\widetilde qv=0,\quad x\in \mathbb{R},\quad 0<t<2\pi,\\
(\Gamma_0v)(t)=\begin{pmatrix}f_1(t) \\ f_2'(t)\end{pmatrix}, \quad 0<t<2\pi,\\
v(\cdot,0)=v_t(\cdot,0)=0.
\end{cases}
\end{equation}
In \cite{MM5} the dynamic IP for (\ref{wave_eqn1}) was studied,
where as a inverse data the authors used the \emph{response
operator}, introduced by the rule
\begin{equation*}
\left(\widetilde
R^TF\right)(t):=\left(\widetilde\Gamma_1v^F\right)(t),\quad t>0.
\end{equation*}
On comparing the representation (\ref{UF_repr}) with one obtained
in \cite{MM5} in Theorem 1, one deduce that for $0<t<2\pi$ the
following equality holds:
\begin{equation}
\label{coinc1}
v^F(x,t)=\begin{cases} u^{F_+}_1(x,t),\quad 0<x<2\pi,\\
u^{F_-}_1(x+2\pi,t),\quad -2\pi<x<0.
\end{cases}
\end{equation}
Moreover, one has that:
\begin{equation}
\label{coinc2} R^{2\pi}F=\Gamma_1u^F=\widetilde
\Gamma_1v^F=\widetilde R^{2\pi}F,\quad 0<t<2\pi.
\end{equation}
Thus we reduced our initial IP to the IP for dynamical system
(\ref{wave_eqn1}) of recovering the potential $\widetilde q(x),$
on the interval $-\pi<x<\pi$ from $\widetilde R^{2\pi}$.

\section{Equations of IP. }

In this section we briefly outline the results of \cite{MM5} in
applying to our situation. Fix a parameter $0<T\leqslant\pi$ and
introduce the \emph{inner space}, the space of states of the
system (\ref{wave_eqn1}) as $\mathcal{H}^T:=L_2(-T,T)$. The
representation (\ref{coinc1}) implies that $v^F(\cdot,T)\in
\mathcal{H}^T$.

A \emph{control operator} $W^T: \mathcal{F}^T\mapsto
\mathcal{H}^T$ is defined by the formula $W^TF:=v^F(\cdot,T)$. The
\emph{reachable set} is defined by the rule
\begin{equation*}
U^T:=W^T\mathcal{F}^T=\left\{v^F(\cdot,T)\,\big|\,  F\in
\mathcal{F}^T\right\}.
\end{equation*}
It will be convenient for us to associate the outer space
$\mathcal{H}^T=L_2(-T,T)$ with a vector space
$L_2(0,T;\mathbb{R}^2)$ by setting for $a\in L_2(-T,T)$ (we keep
the same notation for a function)
\begin{equation*}
a=\begin{pmatrix}a_1(x) \\ a_2(x)\end{pmatrix}\in
L_2(0,T;\mathbb{R}^2),\quad a_1(x):=a(x),\, a_2(x):=a(-x),\, x\in
(0,T).
\end{equation*}
\begin{theorem}
\label{ControlTheor} The control operator is a boundedly
invertible isomorphism between $\mathcal{F}^T$ and
$\mathcal{H}^T$, and $U^T=\mathcal{H}^T$.
\end{theorem}

The \emph{connecting operator} $C^T:\mathcal{F}^T\mapsto
\mathcal{F}^T$ is introduced via the quadratic form:
\begin{equation*}
\left(C^T
F_1,F_2\right)_{\mathcal{F}^T}=\left(v^{F_1}(\cdot,T),v^{F_2}(\cdot,T)\right)_{\mathcal{H}^T}.
\end{equation*}
The crucial fact in the Boundary Control method is that the
connecting operator is expressed in terms of inverse dynamic data:
\begin{theorem}
The connecting operator $C^T$ admits the following representation:
\begin{equation*}
\left(C^TF\right)(t)=\frac{1}{2}\begin{pmatrix}f_1(t)\\
f_2(t)\end{pmatrix}+\int_0^TC(t,s)\begin{pmatrix}
f_1(s)\\f_2(s)\end{pmatrix}\,ds,
\end{equation*}
where
\begin{align*}
&C_{1,1}(t,s)= p_1(2T-t-s)-p_1(|t-s|),\quad p_1(s)=\int_0^s
r_{11}(\alpha)\,d\alpha,\\
&C_{1,2}(t,s)=\widetilde p_1(2T-t-s)-\widetilde p_1(t-s),\quad
\widetilde p_1(s)=\left\{\begin{array}l \int_0^s
r_{12}(\alpha)\,d\alpha,\, s>0,\\
-\int_0^{-s} r_{12}(\alpha)\,d\alpha,\, s<0,
\end{array}
\right.\\
&C_{2,1}(t,s)=-\widetilde r_{21}(t-s)-\widetilde
r_{21}(2T-t-s),\quad \widetilde r_{21}(s)=\left\{\begin{array}l
r_{21}(s),\, s>0,\\
-r_{21}(-s),\, s<0,
\end{array}
\right.\\
&C_{2,2}(t,s)=-r_{22}(|t-s|)-r_{22}(2T-t-s).
\end{align*}
\end{theorem}

\subsection{Krein equations.}

Let $y(x)$ be a solution to the following Cauchy problem:
\begin{equation}
\label{Cauchy_pr} \left\{
\begin{array}l
-y''+\widetilde qy=0,\quad x\in (-T,T),\\
y(0)=0,\,y'(0)=1.
\end{array}
\right.
\end{equation}

We set up the \emph{special control problem}: to find $F\in
\mathcal{F}^T$ such that $W^TF=y$ in $\mathcal{H}^T$. By the
Theorem \ref{ControlTheor}, such a control $F$ exists, but we can
say even more:
\begin{theorem}
The solution to a special control problem is a unique solution to
the following Krein equation:
\begin{equation}
\label{Krein_eqn}
\left(C^TF\right)(t)=(T-t)\begin{pmatrix} 1\\
0\end{pmatrix},\quad t\in (0,T).
\end{equation}
\end{theorem}
Representation formulas (\ref{UF_repr}) and (\ref{coinc1})  imply
that that the solution $F$ to a special control problem satisfies
relations:
\begin{eqnarray*}
y(T)=v^F(T,T)=\frac{1}{2}f_1(0)-\frac{1}{2}f_2(0),\\
y(-T)=v^F(-T,T)=-\frac{1}{2}f_1(0)-\frac{1}{2}f_2(0).
\end{eqnarray*}
Thus solving (\ref{Krein_eqn}) for all $T\in (0,\pi)$, we recover
the solution $y(x)$ to (\ref{Cauchy_pr}) on the interval
$(-\pi,\pi)$. Then the potential $\widetilde q(x),$ $x\in
(-\pi,\pi)$ can be recovered as $\widetilde
q(x)=\frac{y''(x)}{y(x)}$, $x\in (-\pi,\pi)$, and consequently
\begin{equation*}
q(x)=\begin{cases} \widetilde q(x),\quad 0<x<\pi,\\
\widetilde q(x-2\pi),\quad \pi<x<2\pi.
\end{cases}
\end{equation*}

\subsection{Gelfand-Levitan equations}

We introduce the notations:
\begin{eqnarray*}
\label{C_def} C^T=\frac{1}{2}(I+C),\quad
(Cf)(t)=2\int_0^TC(t,s)\begin{pmatrix}
f(s)\\g(s)\end{pmatrix}\,ds,\\
J^T:\mathcal{F}^T\mapsto \mathcal{F}^T,\quad
\left(J^TF\right)(t)=F(T-t),
\end{eqnarray*}
\begin{equation}
\label{C_wid_def} \widetilde C=J^TCJ^T,\quad \left(\widetilde
CF\right)(t)=\int_0^T \widetilde C(t,s)F(s)\,ds.
\end{equation}
Let $m(x,t)\in C^T\left((0,\pi)^2,R^{2\times 2}\right)$ denote a
matrix-valued function such that $m(x,t)=0$ when $x>t.$ In
\cite{MMS} it was proved the following
\begin{theorem} The unique solution to the
Gelfand-Levitan equation
\begin{equation*}
m(x,s)+\widetilde C(x,s)+\int_0^{\pi}\widetilde
C(x,\alpha)m(\alpha,s)\,d\alpha=0,\quad 0<x<s<\pi.
\end{equation*}
where the kernel $\widetilde C$ is defined by (\ref{C_def}),
(\ref{C_wid_def}), determines the potential by the formula:
\begin{eqnarray*}
q(x)=\begin{cases}2\frac{d}{dx}\left(m_{11}(x,x)-m_{12}(x,x)\right),\quad x\in (0,\pi),\\
-2\frac{d}{dx}\left(m_{11}(2\pi-x,2\pi-x)+m_{12}(2\pi-x,2\pi-x)\right),\quad
x\in (\pi,2\pi).
\end{cases}
\end{eqnarray*}
\end{theorem}

\section{Relationship between dynamic and spectral inverse data. }

The problem of finding relationships between different types of
inverse data is very important in inverse problems theory. We can
mention \cite{B01JII,B03,AMR,MM3,MM4,MMS} on some recent results
in this direction. Below we show the relationships between the
dynamic response function, matrix spectral measure and Weyl
matrix.

\subsection{Response function and spectral measure.}
Consider two solution to the equation
\begin{equation}
\label{Eq_Sch} -\phi''+q(x)\phi=\lambda\phi,\quad 0<x<2\pi,
\end{equation}
satisfying the Cauchy data:
\begin{equation*}
\varphi(0,\lambda)=0,\, \varphi'(0,\lambda)=1,\,
\theta(0,\lambda)=1,\,\theta'(0,\lambda)=0.
\end{equation*}
The eigenvalues and normalized eigenfunctions of (\ref{Eq_Sch})
with periodic boundary conditions
\begin{equation}
\label{Eq_Sch1} \phi(0)=\phi(2\pi),\quad \phi'(0)=\phi'(2\pi).
\end{equation}
are denoted by $\{\lambda_n,y_n\}_{n=1}^\infty.$ Let
$\beta_n,\gamma_n\in \mathbb{R}$ be such that
\begin{equation*}
y_n(x)=\beta_n\varphi(x,\lambda_n)-\gamma_n\theta(x,\lambda_n),
\end{equation*}
we point out that there can be eigenvalues of multiplicity two.

We evaluate:
\begin{eqnarray*} y_n(0)=-\gamma_n,\quad
y_n(2\pi)=\beta_n\varphi(2\pi,\lambda_n)+\gamma_n\theta(2\pi,\lambda_n),\\
y_n'(0)=\beta_n,\quad
y_n'(2\pi)=\beta_n\varphi'(2\pi,\lambda_n)+\gamma_n\theta'(2\pi,\lambda_n).
\end{eqnarray*}
Then
\begin{equation}
\Gamma_1y_n=\frac{1}{2}\begin{pmatrix}y_n'(0)+y'_n(2\pi)\\
-y_n(0)-y_n(2\pi)=\end{pmatrix}=\begin{pmatrix}\beta_n\\
\gamma_n\end{pmatrix}.
\end{equation}

Let $F\in \mathcal{F}^T\cap C^\infty_0(0,T;\mathbb{R}^2)$, and
$u^F$ be a solution to (\ref{wave_eqn}). On multiplying
(\ref{wave_eqn}) by $y_n$ and integrating by parts, we get the
following relation:
\begin{eqnarray*}
0=\int_{0}^{2\pi}
u^F_{tt}y_n\,dx-\int_{0}^{2\pi}u^F_{xx}y_n\,dx+\int_{0}^{2\pi}
q(x)u^Fy_n\,dx=\int_{0}^{2\pi} u^F_{tt}y_n\,dx\\
+\left(u^F,Hy_n\right)+\left(\Gamma_1u^F,\Gamma_0y_n\right)_B-\left(\Gamma_0u^F,\Gamma_1y_n\right)_B\\
=\int_{0}^{2\pi}
u^F_{tt}y_n\,dx+\lambda_n\left(u^F,y_n\right)-\left(\begin{pmatrix}f_1(t)\\ f_2'(t)\end{pmatrix},\begin{pmatrix}\beta_n\\
\gamma_n\end{pmatrix}\right)_B.
\end{eqnarray*}
Looking for the solution to (\ref{wave_eqn}) in a form
\begin{equation}
\label{spectr_repr} u^F=\sum_{k=1}^\infty c_k(t)y_k(x),
\end{equation}
we plug (\ref{spectr_repr}) into (\ref{wave_eqn}) and multiply by
$y_n$ to get:
\begin{equation*}
\int_{-N}^N
\sum_{k=1}^\infty c_k''(t)y_k(x)y_n(x)\,dx+\int_{-N}^N\sum_{k=1}^\infty c_k(t)y_k(x)\lambda_ny_n(x)\,dx=\left(\begin{pmatrix}f_1(t)\\ f_2'(t)\end{pmatrix},\begin{pmatrix}\beta_n\\
\gamma_n\end{pmatrix}\right)_B.
\end{equation*}
Thus we obtain that $c_n(t)$, $n\geqslant 1,$ satisfies the
following Cauchy problem:
\begin{equation*}
\left\{ \begin{array}l c_n''(t)+\lambda_nc_n(t)=\left(\begin{pmatrix}f_1(t)\\ f_2'(t)\end{pmatrix},\begin{pmatrix}\beta_n\\
\gamma_n\end{pmatrix}\right)_B,\\
c_n(0)=0,\, c_n'(0)=0.
\end{array}
\right.
\end{equation*}
the solution of which is given by the formula
\begin{equation*}
c_n(t)=\int_0^t\frac{\sin{\sqrt{\lambda_n}(t-s)}}{\sqrt{\lambda_n}}\left(f_1(s)\beta_n+f_2'(s)\gamma_n\right)\,ds.
\end{equation*}
Then for $u^F$ (\ref{spectr_repr}) we have the expansion:
\begin{eqnarray}
u^F(x,t)=\sum_{k=1}^\infty\int_0^t\frac{\sin{\sqrt{\lambda_n}(t-s)}}{\sqrt{\lambda_n}}\left(f_1(s)\beta_n+f_2'(s)\gamma_n\right)\,ds\left(\beta_n\varphi(x,\lambda_n)-\gamma_n\theta(x,\lambda_n)\right)\notag\\
=\sum_{k=1}^\infty\int_0^t\frac{\sin{\sqrt{\lambda_n}(t-s)}}{\sqrt{\lambda_n}}\left(\begin{pmatrix}\beta_n\\
\gamma_n\end{pmatrix}\otimes \begin{pmatrix}\beta_n\\
\gamma_n\end{pmatrix}\begin{pmatrix}f_1(s)
\\ f_2'(s)\end{pmatrix},\begin{pmatrix}\varphi(x,\lambda_n)\\
-\theta(x,\lambda_n)\end{pmatrix}\right)\notag\\
=\int_{-\infty}^\infty\int_0^t\frac{\sin{\sqrt{\lambda}(t-s)}}{\sqrt{\lambda}}\left(d\Sigma(\lambda)\begin{pmatrix}f_1(s)
\\ f_2'(s)\end{pmatrix},\begin{pmatrix}\varphi(x,\lambda)\\
-\theta(x,\lambda)\end{pmatrix}\right).\label{U_F_repr}
\end{eqnarray}
Where $d\Sigma(\lambda)$ is a matrix measure (see \cite{Le})
introduced by the rule:
\begin{equation}
\Sigma(\lambda)=\sum_{\{k\,|\,\lambda_k<\lambda\}}\begin{pmatrix}\beta_n\\
\gamma_n\end{pmatrix}\otimes \begin{pmatrix}\beta_n\\
\gamma_n\end{pmatrix}.
\end{equation}

Thus the response operator $R^T$ is given by
\begin{eqnarray}
(RF)(t)=\Gamma_1v^F=\sum_{k=1}^\infty c_k(t)\Gamma_1y_k=\sum
c_k(t)\begin{pmatrix} \beta_k\\
\gamma_k\end{pmatrix}\label{R_T_reprSP}\\
=\sum_{k=1}^\infty\int_0^t\frac{\sin{\sqrt{\lambda_k}(t-s)}}{\sqrt{\lambda_k}}\left(f_1(s)\beta_k+f_2'(s)\gamma_k\right)\,ds\begin{pmatrix} \beta_k\\
\gamma_k\end{pmatrix}\notag\\
=\int_{-\infty}^\infty\int_0^t
\frac{\sin{\sqrt{\lambda}(t-s)}}{\sqrt{\lambda}}d\Sigma(\lambda)\begin{pmatrix}
f_1(s)\\ f_2'(s)\end{pmatrix}\,ds,\quad 0<t.\notag
\end{eqnarray}

\subsection{Weyl function and response function.}

Let $N_\lambda:=\operatorname{ker}\left(H^*-\lambda I\right)$, we
observe that any $\psi(x,\lambda)\in N_\lambda$ is given by
\begin{equation}
\psi(x,\lambda)=c_1\varphi(x,\lambda)+c_2\theta(x,\lambda).
\end{equation}
We evaluate:
\begin{eqnarray*}
\psi_0=c_2,\quad \psi_{2\pi}=c_1\varphi(2\pi)+c_2\theta(2\pi),\\
\psi'_0=c_1,\quad \psi'_{2\pi}=c_1\varphi'(2\pi)+c_2\theta'(2\pi).
\end{eqnarray*}
Thus the following relations hold:
\begin{eqnarray*}
\Gamma_0\psi=\begin{pmatrix} -\varphi(2\pi) & 1-\theta(2\pi)\\
1-\varphi'(2\pi) & -\theta'(2\pi)\end{pmatrix}\begin{pmatrix}c_1
\\ c_2\end{pmatrix},\\
\Gamma_1\psi= \frac{1}{2}\begin{pmatrix} 1+\varphi'(2\pi) & \theta'(2\pi)\\
-\varphi(2\pi) &
-\left(1+\theta(2\pi)\right)\end{pmatrix}\begin{pmatrix}c_1
\\ c_2\end{pmatrix}.
\end{eqnarray*}
The Weyl matrix is given by (see \cite{BMN})
\begin{equation*}
M(\lambda)=\Gamma_1\left(\Gamma_0|_{N_\lambda}\right)^{-1},
\end{equation*}
so we have:
\begin{equation*}
M(\lambda)=\frac{1}{2}\begin{pmatrix} 1+\varphi'(2\pi) & \theta'(2\pi)\\
-\varphi(2\pi) &
-\left(1+\theta(2\pi)\right)\end{pmatrix}\frac{1}{\det{\Gamma_0}}\begin{pmatrix} -\theta'(2\pi) & -\left(1-\varphi'(2\pi)\right)\\
-\left(1-\theta(2\pi)\right) & -\varphi(2\pi)\end{pmatrix}^T.
\end{equation*}
Evaluating the last expression we get the following formula for
the Weyl matrix:
\begin{equation*}
M(\lambda)=\frac{1}{2\left(F(2\pi,\lambda)-2\right)}\begin{pmatrix} -2\theta'(2\pi,\lambda)\left(1-\varphi'(2\pi,\lambda)\right) & -\varphi'(2\pi,\lambda)+\theta(2\pi,\lambda)\\
-\varphi'(2\pi,\lambda)+\theta(2\pi,\lambda) &
2\varphi(2\pi,\lambda)\end{pmatrix},
\end{equation*}
where
\begin{equation*}
F(x,\lambda)=\varphi'(x,\lambda)+\theta(x,\lambda)
\end{equation*}
is a Lyapunov function.

In \cite{AMR} the authors established the relationship between the
Weyl function and the kernel of dynamic response operator (see
also \cite{MM3,MM4,MMS}). Note that one needs to know the response
for all $t>0$. Then, cf. (\ref{Resp_repr1}):
\begin{equation*}
M(k^2)=\int_0^\infty\left(-\frac{1}{2}\sum_{k=1}^\infty\begin{pmatrix}
\delta'(t-2k\pi)&0\\0&-\delta(t-2k\pi)\end{pmatrix}+\widetilde
R(t)\right)e^{ikt}\,dt,
\end{equation*}
where this equality is understood in a weak sense.

\noindent{\bf Acknowledgments}

The research of Victor Mikhaylov was supported by RFBR
17-01-00529. Alexandr Mikhaylov was supported by RFBR 17-01-00099;
A. S. Mikhaylov and V. S. Mikhaylov were partly supported by RFBR
18-01-00269 and by the Ministry of Education and Science of
Republic of Kazakhstan under grant AP05136197.

\end{document}